\documentclass[12pt]{amsart}
\usepackage{amsmath,amssymb,amsthm,amscd}
\def\qed{\vbox{\hrule
\hbox{\vrule\hbox to 5pt{\vbox to 8pt{\vfil}\hfil}\vrule}\hrule}}

\newcommand{\beg}{\begin{eqnarray*}}
\newcommand{\begn}{\begin{eqnarray}}
\newcommand{\en}{\end{eqnarray*}}
\newcommand{\enn}{\end{eqnarray}}

\newtheorem{prop}{Proposition}[section]
\newtheorem{theo}[prop]{Theorem}
\newtheorem{lemm}[prop]{Lemma}
\newtheorem{coro}[prop]{Corollary}
\newtheorem{rema}[prop]{Remark}

\begin{document}
\title{ Maximal Slice in Anti-de Sitter Space }

\author{ZhenYang Li}
\address{Key Laboratory of Pure and Applied mathematics, School of Mathematics Science, Peking University,
Beijing, 100871, P.R. China.} \email{lzymath@163.com}

 \author{ YuGuang shi}
\address{Key Laboratory of Pure and Applied mathematics, School of Mathematics Science, Peking University,
Beijing, 100871, P.R. China.} \email{ygshi@math.pku.edu.cn}

\renewcommand{\subjclassname}{%
  \textup{2000} Mathematics Subject Classification}

\thanks{The research of  the second author is partially supported by NSF grant of China  and
Fok YingTong Education Foundation}

\begin{abstract}
 In this paper, we prove the existence of maximal slices in anti-de Sitter
 spaces (ADS spaces) with small boundary data at spatial infinity.
 The main argument  is implicit function theorem. We
 also get a necessary and sufficient condition for
  boundary behavior of totally geodesic slice in ADS space.  Moreover, we show that any   isometric and
 maximal embedding of hyperbolic spaces into ADS space
  must be totally geodesic. Together with
 this, we see that most of  maximal slices we get in this paper are not isometric to hyperbolic
 spaces, which implies that the Bernstein Theorem in ADS space fails.
\end{abstract}

\maketitle \markboth{Zhenyang Li and Yuguang Shi } {Maximal slice
in anti-de Sitter space }

\section{Introduction }
\setcounter{equation}{0}
\hspace{0.4cm}
 Finding a minimal surface with the  given boundary  data is an interesting problem in
 Riemannian geometry. Particularly, the existence and
 regularity of the minimal hypersurfaces with a prescribing asymptotic boundary at infinity in hyperbolic space
 $\mathbf{H}^n$ have been discussed in
 \cite{An1}, \cite{An2}, \cite{HL}, \cite{L}, etc.
On the other hand, we know that a
  maximal slice, which is a spacelike hypersurface of a
  Lorentzian manifold and critical point of the induced
  area functional, plays an important role in
  General Relativity. It was used in the first proof of the
  positive mass theorem (\cite{SY}) and in the analysis of the Cauchy problem
  for asymptotically flat spacetimes. Many interesting results for the
  existence of compact maximal slice had been obtained in
  e.g. \cite{B}, \cite{BS}, \cite{G}. For complete noncompact
  case, we have known that there are entire solutions in
  asymptotically flat spacetime (see \cite{B}). It should be pointed  out that a
  complete maximal hypersurface in the Minkowski space must be
  totally geodesic, i.e. a hyperplane( see \cite{CY}). Anti-de Sitter(ADS) space is
  a Lorentzian manifold
  with negative constant sectional curvature, it plays the similar role in
  Lorentzian geometry as hyperbolic space does in Riemannian
  geometry. So, it is nature to study maximal slices in ADS spaces,
  this is one hand; on the other hand, we note that all the
  time slices ( level sets of the time function ) are isometric to
  the $\mathbf{H}^n$ and are totally geodesic, hence are maximal. It
  may be of some interest in view of geometry to find some maximal
  slices which are not totally geodesic. By assuming a global barrier condition in
  asymptotically ADS space, K. Akutagawa proved
  the existence for the entire maximal slice  with certain decay of height
  function at infinity in  \cite{Ak}; in ADS space case, he has also shown
  that, if the height function of the maximal slice  satisfying
  this decay condition at spatial infinity, the maximal slice must be
  time slice (Proposition 3 in  \cite{Ak}).

  In this paper, we obtain
  some maximal slices by implicit function theorem, which can be
  regarded as perturbations of time slices. These maximal slices are $C^{1,1}$ up to the
 boundary. We
 also get a necessary and sufficient condition for
  boundary behavior of totally geodesic slice in ADS space.  Moreover, we show that any   isometric and
 maximal embedding of hyperbolic spaces into ADS space
  must be totally geodesic. Together with
 this, we see that most of  maximal slices we get in this paper are not isometric to hyperbolic
 spaces, which implies that the Bernstein Theorem in ADS space fails.

Indeed, a maximal slice in ADS space satisfies  a second order PDE
in $\mathbf{H^{n}}$(see (1)). Therefore,  it nature to consider the
Dirichlet problem for the maximal slice of  ADS space with infinity
boundary value on $\mathbf{H^{n}}$. We shall address this problem in
the forthcoming papers.

 This paper is organized as follows: In
 Section 2, we derive the equation satisfied by the maximal slices and its corresponding linearized equation.
 In Section 3, we show that the linearized operator is an
 isomorphism between some weighted
H\"older spaces. Hence, using the implicit function theorem, we
prove our main result Theorem 3.1. In Section 4, we prove a
necessary and sufficient condition for the boundary behavior of
totally geodesic slice in ADS space,  we also show that isometric
and maximal embedding of $\mathbf{H}^n$ into ADS spaces is totally
geodesic, by these facts we see that most of  our solutions are not
totally geodesic.

\section{Maximal slice equation in anti-de Sitter
space }
\setcounter{equation}{0}
In this section, we will derive
the maximal slice equation in anti-de Sitter
 space. Let us begin with  some basic facts.
 Suppose  $\mathbf{H}^{n}=(\mathbf{R}^{+}\times \mathbf {S}^{n-1}, \sinh^{-2}\rho
(d\rho^{2}+ d\sigma_{0}^2))$, where $d\sigma_{0}^2$ is the
standard metric on  ${\mathbf S}^{n-1}$. Then $n+1$ dimensional
anti-de Sitter space $V$ can be expressed as a warp product of
$\mathbf{R}$ and $\mathbf{H}^{n}$, namely,
$V=(\mathbf{R}\times\mathbf{H}^{n}, ds^2)$, here
$ds^{2}=-\coth^{2}\rho dt^{2}+\sinh^{-2}\rho (d\rho^{2}+
d\sigma_{0}^2)$. As well known, $V$ is a vacuum solution of
Einstein fields equations with negative cosmological constant. We
denote the canonical connection  in $V$ by $\overline{\nabla}$.
Let $M^{n}$ be a smooth spacelike hypersurface in $V$. The height
function $u\in C^{\infty}(M)$ of $M$ is the restriction of the
time function $t$ to $M$, then $M$ can be regarded  as a graph
over $\mathbf{H}^{n}$, in the following, we assume
$M=\{(x,u(x))|x\in \mathbf{H}^{n}\}$,  and
 $u$ is defined on the whole $V$ by requiring $\frac{\partial}{\partial t}u=0$.

Note that  $M$ is then a level set of $ f(t,x)=t-u(x)$, by a
direct computation, we see that  the future-directed unit normal
vector $N$ to $M$ is
$$N=|\overline{\nabla} f |^{-1}[\overline{\nabla}f ]=\frac{1}{\sqrt{1-\coth^{2}\rho|\nabla u|^{2}}}
(\coth\rho\nabla u+\tanh \rho \frac{\partial }{\partial t}),$$
here and in the sequel, $\nabla$, $div$, and $\Delta$ are the
gradient, divergence and Laplacian operator on $\mathbf{H}^n$
respectively.

Let $\Box$ be the wave operator in $V$,  $H_M$ be the mean curvature
of $M$ in $V$ with respect to $N$, then, by a direct computation, we
see that

$$\Box f =-\Delta u + \tanh \rho \frac{\partial u}{\partial \rho}.$$
On the other hand, we also have

$$\Box f =- NN f -H_M \cdot Nf, $$
thus, we see that

$$
H_{M} =\tanh\rho ~div(\frac{\coth^{2}\rho\nabla u
}{\sqrt{1-\coth^{2}\rho|\nabla u|^{2}}})
$$
If $M$ is maximum, we have

$$
\tanh\rho ~div(\frac{\coth^{2}\rho\nabla u
}{\sqrt{1-\coth^{2}\rho|\nabla u|^{2}}})=0.$$
 One can easily verify that the above equation is equivalent to
\begin{eqnarray}
 div(\frac{\nabla u
}{\sqrt{1-\coth^{2}\rho|\nabla
u|^{2}}})-\frac{2\tanh\rho\frac{\partial u}{\partial \rho}
}{\sqrt{1-\coth^{2}\rho|\nabla u|^{2}}}=0.\end{eqnarray}

By the structure of the equation, we find  that if $u_{\epsilon}$
is the solution of the following equation

\begin{eqnarray} div(\frac{\nabla u
}{\sqrt{1-\epsilon\coth^{2}\rho|\nabla
u|^{2}}})-\frac{2\tanh\rho\frac{\partial u}{\partial \rho}
}{\sqrt{1-\epsilon\coth^{2}\rho|\nabla u|^{2}}}=0,\end{eqnarray}
for some $\epsilon>0$, then $\sqrt{\epsilon}u_{\epsilon}$ is the
solution for equation (1).

In the following, we consider a family of operators:
\begin{equation}
\begin{split}
F(u, \epsilon)&:= div(\frac{\nabla u
}{\sqrt{1-\epsilon\coth^{2}\rho|\nabla
u|^{2}}})-\frac{2\tanh\rho\frac{\partial u}{\partial \rho}
}{\sqrt{1-\epsilon\coth^{2}\rho|\nabla u|^{2}}}\\
&=0,
\end{split}\nonumber
\end{equation}
 and it is easy to see that
$$F(u,0)=\triangle u
-2\tanh\rho \frac{\partial u}{\partial \rho}$$
is the linearize
equation of (1) at its trivial solution $u=0$.

For the purpose of further  discussion, we need to consider the
following Dirichlet problem
\begin{eqnarray}
\left \{
         \begin{array}{lll}
            \mathbf{L}(u):=  \triangle u -2\tanh\rho \frac{\partial u}{\partial \rho}=0,
             \qquad \text{in}\quad \mathbf{H}^n\\
              u |_{\mathbf{S^{n-1}}}=\varphi

        \end{array}  \right.
\end{eqnarray}
where $\varphi$ is a smooth function defined on the infinity
boundary of $\mathbf{H}^{n}$. In (3) and the sequel,
$\mathbf{S^{n-1}}$ is regarded as the infinity boundary of
$\mathbf{H}^n$.
 Besides above facts,  we need to  introduce
the ball model for $\mathbf{H^{n}}$ which is denoted by
$(\mathbf{D^{n}},dS^{2}) $, here, $\mathbf{D^{n}}$ is the unit
ball in $\mathbf{R^n}$, and $dS^2$ is the standard hyperbolic
metric which is defined as following:

$$dS^{2}=\tau^{-2}\displaystyle\sum_{i=1}^{n}(dx^{i})^{2},$$
where $\tau(x)=\frac{1}{2}(1-|x|^{2})$ and
$\displaystyle\sum_{i=1}^{n}(dx^{i})^{2}$ is the Euclidean metric.
The relation between  $\rho$ and $\tau$ can be expressed as
$\sinh\rho=\frac{\tau}{r}$, where $r(x)=|x|$ is the Euclidean
distance from the origin. Hence the equation (1), (2) and  (3) can
be written as
\begin{eqnarray}
 div(\frac{\nabla u
}{\sqrt{1-(\frac{1-\tau}{\tau})^{2}|\nabla
u|^{2}}})+\frac{2\tau}{1-\tau}\frac{\displaystyle\sum_{i=1}^{n}x^i
\frac{\partial u}{\partial x^i}
}{\sqrt{1-(\frac{1-\tau}{\tau})^{2}|\nabla
u|^{2}}}=0,\nonumber\end{eqnarray}

\begin{eqnarray}
 div(\frac{\nabla u
}{\sqrt{1-\epsilon(\frac{1-\tau}{\tau})^{2}|\nabla
u|^{2}}})+\frac{2\tau}{1-\tau}\frac{\displaystyle\sum_{i=1}^{n}x^i
\frac{\partial u}{\partial x^i}
}{\sqrt{1-\epsilon(\frac{1-\tau}{\tau})^{2}|\nabla
u|^{2}}}=0.\nonumber\end{eqnarray} and
\begin{eqnarray}
\left \{
         \begin{array}{lll}
            \mathbf{L}(u):= \triangle u + \frac{2\tau}{1-\tau}
\displaystyle\sum_{i=1}^{n}x^i \frac{\partial u}{\partial x^i}=0 \qquad \text{in}\quad \mathbf{H}^n\\
              u |_{\mathbf{S^{n-1}}}=\varphi

        \end{array}  \right.
\end{eqnarray}
respectively.

\section{Existence of maximal slice, Weighted H\"older space and analysis of linearize  equation}

In this section, we will prove the existence of maximal slices in
$V$ with certain boundary data at infinity. More specifically, we
are going to show

\begin{theo}
For any $\varphi \in C^{4,\alpha}(\mathbf {S}^{n-1})$, there is a
$\delta=\delta(\varphi)>0$ so that for any $\epsilon \in (0,
\delta)$, the following Dirichlet problem
\begin{eqnarray}
\left \{
         \begin{array}{lll}
            div(\frac{\nabla u
}{\sqrt{1-(\frac{1-\tau}{\tau})^{2}|\nabla
u|^{2}}})+\frac{2\tau}{1-\tau}\frac{\displaystyle\sum_{i=1}^{n}x^i
\frac{\partial u}{\partial x^i}
}{\sqrt{1-(\frac{1-\tau}{\tau})^{2}|\nabla
u|^{2}}}=0, \qquad\text{in} \quad \mathbf{H}^n\\
              u |_{\mathbf{S^{n-1}}}=\sqrt{\epsilon}\varphi

        \end{array}  \right.
\end{eqnarray}
admits a solution $u \in C^2 (\mathbf{D}^n)$ with $\|u\|_{C^2
(\mathbf{D}^n)}\leq C$, here $C$ is a constant depends on $\varphi$.
\end{theo}

\begin{rema}
\begin{enumerate}
\item In Theorem 3.1, we adopt the ball model for $\mathbf{H}^n$,
and $u$ is regarded as a function defined on $\mathbf{D}^n$.

  \item  The second fundamental form of the solution we get in Theorem 3.1
decays as $O(\tau^{2})$ as $\tau$ goes to $0$. And we conjecture
that the solution with  the second fundamental form  faster than
quadratic decay must be $\mathbf{H}^n$.
\end{enumerate}

\end{rema}

In order to prove Theorem 3.1, we will get some basic estimates of
the linear equation by which we are able to show the corresponding
linear elliptic operator is a linear isomorphism between some
function spaces. To do this, let us first introduce a kind of
weighted H\"older spaces (for more details, please see \cite{GL}.).
In the following, we will define weighted H\"older spaces on
$\Omega\subset\mathbf{D^{n}}$,  for $0\leq k\in \mathbf{Z}$ let
$C^{k}(\overline{\Omega}) $  be the usual Banach spaces of $k$ times
continuously differential functions on $\overline{\Omega}$, and
$0<\alpha<1$ denote by $C^{k,\alpha}(\overline{\Omega}) $ the
subspace of functions whose $k-$th derivatives satisfy a uniform
H\"older condition of order $\alpha$, with the usual norms denoted
by $||\cdot||_{k;\Omega}$, $||\cdot||_{k,\alpha;\Omega}$
respectively. And denote $C^{k}(\Omega)$ and $C^{k,\alpha}(\Omega)$
the linear space of functions satisfying the corresponding estimates
uniformly on compact subsets of $\Omega$. For $s\in \mathbf{R}$
define

\begin{eqnarray}
||w||_{k,0;\Omega}^{(s)}=\displaystyle\sum_{l=0}^{k}\displaystyle\sum_{|\gamma|=l}\|\tau^{-s+l}
\partial^{\gamma}w\|_{L^{\infty}(\Omega)},\nonumber
\end{eqnarray}
where for any multi-index $\gamma$,\hspace*{0.15cm}$
\partial^{\gamma}=\frac{\partial^{|\gamma|}}{\partial
x^{\gamma}}$;\hspace*{0.15cm} and for $0<\alpha<1$ define

\begin{equation}
\begin{split}
||w||_{k,\alpha;\Omega}^{(s)}&=||w||_{k,0;\Omega}^{(s)}\\
&+\displaystyle\sum_{|\gamma|=k}\displaystyle\sup_{x,y\in
\Omega}[\min(\tau^{-s+k+\alpha}(x),\tau^{-s+k+\alpha}(y))
\frac{|\partial^{\gamma}w(x)-\partial^{\gamma}w(y)|}{|x-y|^{\alpha}}]\nonumber
\end{split}\nonumber
\end{equation}
Let $\Lambda^{s}_{k,\alpha;\Omega}=\{w\in
C^{k,\alpha}(\Omega)\mid\hspace*{0.15cm}
||w||^{(s)}_{{k,\alpha;\Omega}}<+\infty\}$, which is a Banach
space. For $x\in\mathbf{H^{n}}$, let $B(x)$ be the open Euclidean
ball with center $x$ and radius $\frac{1}{3}\tau(x)$. It is clear
that

\begin{lemm} For any $\Omega' \subset \Omega \subset \mathbf{H^{n}}$, we have
$ \Lambda^{s}_{k,\alpha; \Omega}\subset\Lambda^{s}_{k,\alpha;
\Omega'}$, and
\begin{eqnarray}
||w||_{k,\alpha;\Omega'}^{(s)}\leq||w||_{k,\alpha;\Omega}^{(s)}\nonumber
\end{eqnarray}
 for any $ w\in\Lambda^{s}_{k,\alpha;
\Omega}$; also, For any $\Omega_m \subset
\Omega_{m+1}\subset\mathbf{H^{n}}$ and $\displaystyle\bigcup_m
\Omega_m = \mathbf{H^{n}}$,we have
\begin{eqnarray}
||w||_{k,\alpha;\mathbf{H^{n}}}^{(s)}\leq
\displaystyle\sup_{m}||w||_{k,\alpha;\Omega_m}^{(s)}\nonumber
\end{eqnarray}
for any $ w\in\Lambda^{s}_{k,\alpha; \mathbf{H}^n}$.
\end{lemm}

\begin{proof}
By the definition, we see that for any $w\in
\Lambda^{s}_{k,\alpha;\Omega}$,

$$\|\tau^{-s+l}\partial^\gamma w \|_{L^\infty (\Omega')}\leq
\|\tau^{-s+l}\partial^\gamma w \|_{L^\infty (\Omega)}, $$ and

\begin{equation}
\begin{split}
&\displaystyle\sup_{x,y\in
\Omega'}[\min(\tau^{-s+k+\alpha}(x),\tau^{-s+k+\alpha}(y))
\frac{|\partial^{\gamma}w(x)-\partial^{\gamma}w(y)|}{|x-y|^{\alpha}}]\\
&\leq \displaystyle\sup_{x,y\in
\Omega}[\min(\tau^{-s+k+\alpha}(x),\tau^{-s+k+\alpha}(y))
\frac{|\partial^{\gamma}w(x)-\partial^{\gamma}w(y)|}{|x-y|^{\alpha}}]
\end{split}\nonumber
\end{equation}
Thus, we see that
\begin{eqnarray}
||w||_{k,\alpha;\Omega'}^{(s)}\leq||w||_{k,\alpha;\Omega}^{(s)};
\nonumber
\end{eqnarray}
and thus $ \Lambda^{s}_{k,\alpha;
\Omega'}\subset\Lambda^{s}_{k,\alpha; \Omega}$. On the other hand,
for any $\epsilon>0$, there is an $x \in \mathbf{H^{n}}$, we may
assume $x \in \Omega_m$ so that

$$|\tau^{-s+l}\partial^\gamma w(x)|> \|\tau^{-s+l}\partial^\gamma w\|_{L^\infty (\mathbf{H^{n}})}-\epsilon,$$
and
$$|\tau^{-s+l}\partial^\gamma w(x)|\leq \|\tau^{-s+l}\partial^\gamma w\|_{L^\infty (\Omega_m)}$$
hence, we see
$$\|\tau^{-s+l}\partial^\gamma w\|_{L^\infty (\mathbf{H^{n}})}-\epsilon\leq
\|\tau^{-s+l}\partial^\gamma w\|_{L^\infty (\Omega_m)}
$$
By the same arguments, we have
\begin{equation}
\begin{split}
&\displaystyle\sup_{x,y\in
\mathbf{H^{n}}}[\min(\tau^{-s+k+\alpha}(x),\tau^{-s+k+\alpha}(y))
\frac{|\partial^{\gamma}w(x)-\partial^{\gamma}w(y)|}{|x-y|^{\alpha}}]\\
&\leq \displaystyle\sup_{x,y\in
\Omega_m}[\min(\tau^{-s+k+\alpha}(x),\tau^{-s+k+\alpha}(y))
\frac{|\partial^{\gamma}w(x)-\partial^{\gamma}w(y)|}{|x-y|^{\alpha}}]+\epsilon
\end{split}\nonumber
\end{equation}

Thus, for any $\epsilon>0$, and sufficiently large $m$, we have
\begin{eqnarray}
||w||_{k,\alpha;\mathbf{H^{n}}}^{(s)}\leq
\displaystyle\sup_{m}||w||_{k,\alpha;\Omega_m}^{(s)}+ \epsilon
\nonumber
\end{eqnarray}
which implies the conclusion is true.
\end{proof}

The following lemma is the same as Lemma 3.1 in \cite{GL}.
\begin{lemm}

For $x\in\Omega$, then
\begin{eqnarray}
||w||_{k,\alpha;B(x)\cap\Omega}^{(s)}\leq||w||_{k,\alpha;\Omega}^{(s)}
\nonumber
\end{eqnarray}
and,
\begin{eqnarray}
||w||_{k,\alpha;\Omega}^{(s)}\leq
C\displaystyle\sup_{x\in\Omega}||w||_{k,\alpha;B(x)\cap\Omega}^{(s)}
\nonumber
\end{eqnarray}
where $C$  depends only on $k$.
\end{lemm}
Let $B(0)$ denote the open Euclidean ball with center $0$ and
radius $\frac{1}{3}$, and for $x\in\mathbf{H^{n}}$ define
$\psi_{x}:B(0)\rightarrow B(x)$ by
\begin{eqnarray}y:=\psi_{x}(z)=x+\tau(x)z.\end{eqnarray}
If $y\in B(x)$, then
\begin{eqnarray}\frac{1}{10}\tau(x)\leq\tau(y)\leq40\tau(x).\end{eqnarray}
Therefore, there exist a universal constant $\Lambda_{1}$ such that

\begin{equation}
\begin{split}
&\Lambda_{1}^{-1}\tau^{-s+l}(x)||\partial^{\gamma}w||_{L^{\infty}(B(x)\cap\Omega)}\leq||\tau^{-s+l}
\partial^{\gamma}w||_{L^{\infty}(B(x)\cap\Omega)}\\
&\leq
\Lambda_{1}\tau^{-s+l}(x)||\partial^{\gamma}w||_{L^{\infty}(B(x)\cap\Omega)}
\end{split}\nonumber
\end{equation}
where $\Lambda_{1}$   depends only  on $s$ and $l$.   Let
$v(z)=w\circ\psi_{x}(z)$, we have
$\frac{\partial^{\gamma}v}{\partial
z^{\gamma}}=\tau^{l}(x)\frac{\partial^{\gamma}w}{\partial
y^{\gamma}}$ for $|\gamma|=l$. So for any $y\in B(x)\cap\Omega$,
we can show
\begin{eqnarray}
\tau^{-s+l}(x)\partial^{\gamma}w(y)=
\tau^{-s}(x)\partial^{\gamma}v(z).\nonumber
\end{eqnarray}
Using (6) and (7),  one can conclude that
\begin{eqnarray}
\begin{array}{ccc}&\Lambda_{1}^{-1}\tau^{-s}(x)||\partial^{\gamma}v||_{L^{\infty}(\psi_{x}^{-1}(B(x)\cap\Omega))}\leq
||\tau^{-s+l}\partial^{\gamma}w||_{L^{\infty}(B(x)\cap\Omega)}&\\&\leq
\Lambda_{1}\tau^{-s}(x)||\partial^{\gamma}v||_{L^{\infty}(\psi_{x}^{-1}(B(x)\cap\Omega))}&\end{array}\nonumber
\end{eqnarray}
By this, it follows that
\begin{eqnarray}
\begin{array}{ccc}&\Lambda^{-1}\tau^{-s}(x)||v||_{k,\alpha;
\psi_{x}^{-1}(B(x)\cap\Omega)}\leq
||w||_{k,\alpha;B(x)\cap\Omega}^{(s)}&\\&\leq
\Lambda\tau^{-s}(x)||v||_{k,\alpha;
\psi_{x}^{-1}(B(x)\cap\Omega)}&
\end{array}\end{eqnarray}
where $\Lambda$ is only depended on $k$ , $\alpha$ and $s$.\\\\

Next, consider the following
\begin{equation}
\left \{
         \begin{array}{lll}
            \mathbf{L}(u)=  \triangle u -2\tanh\rho \frac{\partial u}{\partial \rho}
            =\triangle u+2\frac{\tau(y)}{1-\tau(y)}  \displaystyle\sum_{i=1}^{n}
 y^{i}\frac{\partial }{\partial y^{i}}u=\eta\\
              u |_{\mathbf{S^{n-1}}}=0

        \end{array}  \right.
\end{equation}
here, $\eta \in \Lambda^{s}_{0,\alpha;\mathbf{H^n}}$ where $s$ is
to be determined later.

\begin{lemm}Suppose $u\in
C^{2}(\mathbf{H^{n}})\bigcap\Lambda^{s}_{0,0;\mathbf{H^{n}}}$ is a
solution for (9) , $\eta\in\Lambda^{s}_{k,\alpha;\mathbf{H^{n}}}$.
Then
\begin{eqnarray}||u||_{k+2, \alpha; \mathbf{H^{n}}}^{(s)}\leq C(||\eta||_{k ,\alpha;
\mathbf{H^{n}}}^{(s)}+||u||_{0,0;\mathbf{H^{n}}}^{(s)})\nonumber\end{eqnarray}
here $C=C(k,\alpha)$.
\end{lemm}

\begin{proof} It is easy to see that (9) is equivalent to
\begin{eqnarray}
\hspace{1cm}\left\{\begin{array}{lll} \tau^{2}(y) \triangle_{0}
u+\tau(y)(n-2+\frac{2}{1-\tau(y)})\displaystyle\sum_{i=1}^{n}
y^{i}\frac{\partial }{\partial y^{i}}u
=\eta& \mbox{in}\hspace{0.2cm} \mathbf{D^{n}}\\
                       u|_{_{\mathbf{S^{n-1}}}}=0, &
\end{array}\right. \end{eqnarray}\\
where $\triangle_{0}$ is the standard Laplacian for
$\mathbf{D^{n}}\subset \mathbf{R^{n}}$. Suppose
$v(z)=u\circ\psi_{x}(z)$ for $ \forall z\in B(0)$, then (10)
becomes

\begin{eqnarray}
\hspace{1cm}\left\{\begin{array}{lll}
\frac{\tau^{2}(y)}{\tau^{2}(x)}
\triangle_{0}v+\frac{\tau(y)}{\tau(x)}(n-2+\frac{2}{1-\tau(y)})
\displaystyle\sum_{i=1}^{n}y^{i}
\frac{\partial }{\partial z^{i}}v=\eta&\mbox{in}\hspace{0.2cm} B(0)  \\
                       u|_{_{\mathbf{S^{n-1}}}}=0. &

\end{array}\right. \end{eqnarray}\\
Let $B'(0)$ and $B'(x)$ denote the open Euclidean balls with
center $0$ and radius $\frac{1}{4}$ and with center $x$ and radius
$\frac{1}{4}\tau(x)$ respectively. Since
$\frac{1}{100}\leq\frac{\tau^{2}(y)}{\tau^{2}(x)}\leq160$ when
$y\in B(x)$ and $n\leq n-2+\frac{2}{1-\tau(y)}\leq n+2$, it
follows that (11) is uniformly elliptic equation on $B(0)$. Hence
by standard Schauder theory, we have
\begin{eqnarray}
||v||_{k+2,\alpha;B'(0)}\leq C
(||\eta\circ\psi_{x}||_{k,\alpha;B(0)}+||v||_{0,0;B(0)})
\nonumber\end{eqnarray}
 where $C$ only depends on $k,\alpha$.
Choose $\Omega'\subset\Omega$ such that $B(x)\subset\Omega$ for
any $x\in\Omega'$. Applying (8) and Lemma 3.4,  we obtain
\begin{eqnarray}
\begin{array}{llll}
||u||_{k+2, \alpha; \Omega'}^{(s)}&\leq&
C\displaystyle\sup_{x\in\Omega'}\tau^{-s}(x)||u\circ\psi_{x}||_{k+2,\alpha;\psi_{x}^{-1}(B'(x)\cap\Omega')}\\
&\leq&C\displaystyle\sup_{x\in\Omega'}\tau^{-s}(x)||u\circ\psi_{x}||_{k+2,\alpha;B'(0)}\\
&\leq&C\displaystyle\sup_{x\in\Omega'}\tau^{-s}(x)(||\eta\circ\psi_{x}||_{k,\alpha;B(0)}+||v||_{0,0;B(0)})\\
&\leq&C(||\eta||_{k ,\alpha;
\Omega}^{(s)}+||u||_{0,0;\Omega}^{(s)})\\
&\leq&C(||\eta||_{k ,\alpha;
\mathbf{H^{n}}}^{(s)}+||u||_{0,0;\mathbf{H^{n}}}^{(s)})
\end{array} \nonumber\end{eqnarray}\\
Therefore lemma follows from Lemma 3.3.\end{proof}.

\begin{prop} For any  $\eta\in\Lambda^{s}_{0,\alpha;\mathbf{H^{n}}}$ , there exists $u\in
C^{2}(\mathbf{H^{n}})\bigcap\Lambda^{s}_{0,0;\mathbf{H^{n}}}$ for
$0\leq s<n+1$ satisfying (9), and \hspace{0.1cm}$
\|u\|^{(s)}_{0,0;\mathbf{H^{n}}}\leq
C\|\eta\|^{(s)}_{0,\alpha;\mathbf{H^{n}}}$ with $C$ depending on
$s$. Moreover, $u\in \Lambda^s _{2, \alpha; \mathbf{H}^n}$ with
$$\|u\|_{2, \alpha; \mathbf{H}^n}^{(s)}\leq C \|\eta\|_{0, \alpha; \mathbf{H}^n}^{(s)},$$
here $C$ is a constant depends only on $s$ and $\alpha$.
\end{prop}

\begin{proof}Let $\{\Omega_{m}\}_{m=1}^{\infty}$ be an exhausting
sequence such that $\Omega_{m}\subset\Omega_{m+1}$ and
$\displaystyle\bigcup_{m}\Omega_{m}=\mathbf{H^{n}}$. Let $w_{m}$
be the solution for the following equation:
\begin{eqnarray}
\left\{\begin{array}{lll}
\mathbf{L}(w_{m})= \eta\hspace{1cm}\mbox{in}\hspace{0.2cm}\Omega_{m}&  \\
w_{m}|_{_{\partial\Omega_{m}}}=0 .&
\end{array}\right. \nonumber\end{eqnarray}\\
Hence $w_{m}\in C^{2,\alpha}(\overline{\Omega}_{m})$ since
$\eta\in C^{0,\alpha}(\overline{\Omega}_{m})$(\cite{GL}).

Set $\phi=\tau^{s}$,
\begin{eqnarray}
\begin{array}{lll}
\mathbf{L}(\phi)&=&-s(2s-n+2)\tau^{s+1}+s(s-n-1)\tau^{s}+\frac{2}{1-\tau}s\tau^{s+1}\\
&\leq&-s(2s-n-2)\tau^{s+1}+s(s-n-1)\tau^{s}
\end{array}\nonumber \end{eqnarray}
since $\frac{2}{1-\tau}\leq4$ and $s>0$.  For $0\leq s<n+1$, it is
easy to check  that
\begin{eqnarray}
\mathbf{L}(\phi)=\triangle\phi-\frac{2\tau(y)(1-2\tau(y))}{1-\tau(y)}
\frac{\partial }{\partial
\tau}\phi\leq-\delta\phi\nonumber\end{eqnarray} for some constant
$\delta>0$ only depended on $s$.
  On the other
hand,  we have $|\eta|\leq C\tau^{s}$ where
$C=\|\eta\|^{(s)}_{0,\alpha;\mathbf{H^{n}}}$ since
$\eta\in\Lambda^{s}_{0,\alpha;\mathbf{H^{n}}}$. We choose  a
constant $C_{1}=C/\delta$ such that
\begin{eqnarray}
\left\{\begin{array}{lll}
\mathbf{L}(w_{m})\geq \mathbf{L}(C_{1}\phi) \hspace{1cm}\mbox{in}\hspace{0.2cm}\Omega_{m}&\\
                       (C_{1}\phi-w_{m})|_{_{\partial\Omega_{m}}}\geq0 &
\end{array}\right. \nonumber\end{eqnarray}\\
By maximum principle, $w_{m}\leq C_{1}\tau^{s}$. By the same
arguments, we may get the lower bound of $w_m$, hence, $|w_m|\leq
C_1 \tau^s$. Therefore $w_{m}$ convergence to a function $u\in
C^{2}(\mathbf{H}^{n})\bigcap\Lambda^{s}_{0,0;\mathbf{H}^{n}}$
which solves (9), and we have $
\|u\|^{(s)}_{0,0;\mathbf{H^{n}}}\leq
C\|\eta\|^{(s)}_{0,\alpha;\mathbf{H^{n}}}$ where $C$ depends only
on $s$. By Lemma 3.5, we know $u\in \Lambda^s _{2, \alpha;
\mathbf{H}^n}$ with
$$\|u\|_{2, \alpha; \mathbf{H}^n}^{(s)} \leq C \|\eta\|_{0, \alpha; \mathbf{H}^n}^{(s)},$$
where $C=C(s,\alpha)$.
\end{proof}

Now by the  Lemma 3.5 and Proposition 3.6, We have the following:

\begin{theo}
The operator
$\mathbf{L}:\Lambda^{s}_{2,\alpha;\mathbf{H^{n}}}\rightarrow\Lambda^{s}_{0,\alpha;\mathbf{H^{n}}}$
defined in (9) is an isomorphism, where $0\leq s<n+1$.\\\\
\end{theo}

\begin{coro} For any $\varphi \in C^{4, \alpha}(\mathbf{S}^{n-1})$,
Dirichlet problem (3)(or (4)) has a solution.
\end{coro}
\begin{proof}
We use the cylindrical coordinate system $(\rho,\theta)$. And
extend $\varphi$ as $\varphi(\rho,\theta)=\varphi(\theta)$, for
$\theta\in\mathbf{S^{n-1}}$ and small $\rho$. Then let
$f(\rho,\theta)\in C^{2,\alpha}(\mathbf{H^{n}})$ such that for
some small $\rho$,
$f(\rho,\theta)=\varphi+\frac{1}{2(n-1)}\rho^{2}\triangle_{\mathbf{S^{n-1}}}\varphi$,
here, $\triangle_{\mathbf{S^{n-1}}}$ is the Laplacian operator on
$\mathbf{S^{n-1}}$. Put $f$ into left side of (3), one can see
that
\begin{equation}
 \begin{array}{lll}
 \mathbf{L}(f)&=&\frac{1}{n-1}\sinh^{2}\rho
 \triangle_{\mathbf{S^{n-1}}}\varphi-\frac{1}{n-1}((n-2)\sinh\rho\cosh\rho+2\tanh\rho)\\
 &&\cdot\rho
 \triangle_{\mathbf{S^{n-1}}}\varphi+\sinh^{2}\rho\triangle_{\mathbf{S^{n-1}}}\varphi
 +\frac{1}{2(n-2)}\rho^{2}\sinh^{2}\rho\triangle^{2}_{\mathbf{S^{n-1}}}\varphi\\
 &=& O(\rho^{4})=O(\tau^{4}) \hspace{1cm}\mbox{as}\hspace{0.2cm}\tau\rightarrow0.\nonumber\end{array}\end{equation}
 Because $\mathbf{L}(f)$ is $C^{0,\alpha}$ in any compact subset and
 behave like $\tau^{4}$ near boundary, we conclude that
 $\mathbf{L}(f)\in\Lambda^{s}_{0,\alpha;\mathbf{H^{n}}}$, for any $s\leq4$.
 Then the corollary follows from the Proposition 3.6.

\end{proof}

Now we are in the position to prove Theorem 3.1. By Corollary 3.8,
(4) has a solution u satisfying $u-f\in
\Lambda^{s}_{2,\alpha;\mathbf{H^{n}}}$ for some $s\in[0,4)$ where
$f$ is given as in the proof of Corollary 3.8. Since (4) is a
linear equation, we can multiply $\varphi$ by a suitable constant
such that the corresponding solution $u$ satisfying
$v=\frac{1}{\sqrt{1-(\frac{1-\tau}{\tau})^{2}|\nabla
u|^{2}}}<+\infty$ , that is, $u$ is spacelike. And define
$$\Xi_{A}=\{w\in \Lambda^{2}_{2,\alpha;\mathbf{H^{n}}}\mid\hspace*{0.15cm}
\frac{1}{\sqrt{1-(\frac{1-\tau}{\tau})^{2}|\nabla (w+u)|^{2}}}<
A<+\infty\}\subset \Lambda^{2}_{2,\alpha;\mathbf{H^{n}}}$$
Obviously $\Xi_{A}$ is nonempty open set of
$\Lambda^{2}_{2,\alpha;\mathbf{H^{n}}}$, since $0\in\Xi_{A}$.
Define an operator

$$H(\cdot,\cdot):(-1,+1)\times \Xi_{A}\rightarrow\Lambda^{2}_{0,\alpha;\mathbf{H^{n}}}$$
by

\begin{eqnarray}
H(\epsilon, w)&:=& div(\frac{\nabla (w+u)
}{\sqrt{1-\epsilon(\frac{1-\tau}{\tau})^{2}|\nabla
(w+u)|^{2}}})\nonumber\\
&+&\frac{2\tau}{1-\tau}\frac{\displaystyle\sum_{i=1}^{n}x^i
\frac{\partial (w+u)}{\partial x^i}
}{\sqrt{1-\epsilon(\frac{1-\tau}{\tau})^{2}|\nabla
(w+u)|^{2}}}\nonumber\\
&=&\frac{\triangle
(w+u)}{\sqrt{1-\epsilon(\frac{1-\tau}{\tau})^{2}|\nabla
(w+u)|^{2}}}\nonumber\\
&+&<\nabla\frac{1}{\sqrt{1-\epsilon(\frac{1-\tau}{\tau})^{2}|\nabla
(w+u)|^{2}}},
\nabla(w+u)>\nonumber\\
&&+\frac{2\tau}{1-\tau}\frac{\displaystyle\sum_{i=1}^{n}x^i
\frac{\partial (w+u)}{\partial x^i}
}{\sqrt{1-\epsilon(\frac{1-\tau}{\tau})^{2}|\nabla
(w+u)|^{2}}}\nonumber
\end{eqnarray}
From Corollary 3.8, we have $H(0,0)=0$. By a direct computation,
we see that $H$ is a smooth operator, and for any
$h\in\Lambda^{2}_{2,\alpha;\mathbf{H^{n}}}$,
$$\frac{\partial}{\partial
t}H(0,th)|_{t=0}=\triangle h+2\frac{\tau(y)}{1-\tau(y)}
\displaystyle\sum_{i=1}^{n} x^{i}\frac{\partial }{\partial
x^{i}}h$$ It follows that the map $\frac{\partial}{\partial
t}H(0,th)|_{t=0}=\mathbf{L}:\Lambda^{2}_{2,\alpha;\mathbf{H^{n}}}\rightarrow\Lambda^{2}_{0,\alpha;\mathbf{H^{n}}}$
is an isomorphism from Theorem 3.7. Now by the implicit  function
theorem(cf. \cite{GT}), we can conclude that (5) has a solution
whose difference by $u$ is in
$\Lambda^{2}_{2,\alpha;\mathbf{H^{n}}}$ and boundary data is given
by small $\sqrt{\epsilon}\varphi$. Thus we finish to prove Theorem
3.1.

\section{Boundary behavior of totally geodesic slice of ADS space}

In this section, we will show that any isometric and maximal
embedding of $\mathbf{H}^n$ into ADS space is totally geodesic,
and moreover, we will give a sufficient and necessary condition of
the boundary value of the height function for totally geodesic
slice. Together with Theorem 3.1, we know that the Bernstein
Theorem in ADS spacetime fails. Let's begin with the following

\begin{prop}
If a hyperbolic space is isometrically immersed in the anti-de
Sitter space as its maximal hypersurface, it must be totally
geodesic.
\end{prop}
\begin{proof}Suppose that $M$ is a hyperbolic space, which is also a maximal hypersurface
of anti-de Sitter space V. We choose a local field of Lorentzian
orthonormal frames $e_{0}, e_{1}, \ldots, e_{n}$ in $V$ such that,
at each point of $M$, $e_{1}, \ldots, e_{n}$ spans the tangent
space of $M$ and $e_{0}$ is the unit timelike normal vector for
$M$. We make use of the following convention on the ranges of
indices:
$$0\leq\alpha, \beta, \gamma,\ldots\leq n,\hspace{0.2cm} 1\leq i,j,k,\ldots\leq n.$$
Let $\omega_{0}, \omega_{1},\ldots,\omega_{n}$ be the dual frame
field. Then the structure equations of $V$ are
\begin{eqnarray}
\left \{\begin{array}{lll}
d\omega_{0}&=&-\omega_{0i}\bigwedge\omega_{i}\\
d\omega_{i}&=&\omega_{i0}\bigwedge\omega_{0}-\omega_{ik}\bigwedge\omega_{k}\hspace{1cm}
\omega_{\alpha\beta}+\omega_{\beta\alpha}=0\\
d\omega_{0i}&=&-\omega_{0k}\bigwedge\omega_{ki}-K_{0i0j}\omega_{0}\bigwedge\omega_{j}+
\frac{1}{2}K_{0ijk}\omega_{j}\bigwedge\omega_{k}\\
d\omega_{ij}&=&\omega_{i0}\bigwedge\omega_{0j}-\omega_{ik}\bigwedge\omega_{kj}-K_{ijok}\omega_{0}\bigwedge\omega_{k}
+\frac{1}{2}K_{ijkl}\omega_{k}\bigwedge\omega_{l}
\end{array}\right.\nonumber\end{eqnarray}
where $K_{\alpha\beta\gamma\delta}$ is the curvature tensor for
$V$. We restrict these forms to $M$,  then $\omega_{0}=0$. We may
put $\omega_{0i}=h_{ij}\omega_{j}$, where $h_{ij}$ is the
components of the second fundamental form of $M$. And we also have
the structure equation for $M$:
\begin{eqnarray}
\left \{\begin{array}{lll}
d\omega_{i}&=&-\omega_{ik}\bigwedge\omega_{k}\hspace{1cm}
\omega_{ij}+\omega_{ji}=0\\
d\omega_{ij}&=&-\omega_{ik}\bigwedge\omega_{kj}
+\frac{1}{2}R_{ijkl}\omega_{k}\bigwedge\omega_{l}
\end{array}\right.\nonumber\end{eqnarray}
where $R_{ijkl}$ is the curvature tensor of $M$. Hence we have the
Gauss formula,
$$R_{ijkl}=K_{ijkl}-h_{ik}h_{jl}+h_{il}h_{jk}.$$
Since both $V$ and $M$ have constant sectional curvature $-1$, we
can see that$$h_{ii}h_{jj}-h_{ij}^{2}=0$$ for $i\neq j$. Using the
fact that $M$ is maximal, i.e.
$\displaystyle\sum_{i=1}^{n}h_{ii}=0 $,  we have
$$0=h_{jj}\displaystyle\sum_{i\neq
j}h_{ii}-\displaystyle\sum_{i\neq
j}h_{ij}^{2}=-h_{jj}^{2}-\displaystyle\sum_{i\neq j}h_{ij}^{2}.$$
 it follows
 that $M$ is totally geodesic.
\end{proof}

Now, we are in the position to study the boundary behavior of
totally geodesic slice of ADS space $V$. For simplicity, we only
consider the case that $dim V=4$.

Let $\mathbf{R}^{5}_{2}$ be $5$ dimensional semi-Euclidean space,
that is, it is a vector space with the inner product
$<X,Y>=x_{1}y_{1}+x_{2}y_{2}+x_{3}y_{3}-x_{4}y_{4}-x_{5}y_{5}$,
where $X=(x_{1},x_{2},x_{3},x_{4},x_{5})$ and
$Y=(y_{1},y_{2},y_{3},y_{4},y_{5})$. Denote the connection in
$\mathbf{R}^{5}_{2}$ by $\widetilde{\nabla}$. It is  well known that
$4$ dimensional anti-de Sitter space $V$ is a totally umbilical
hypersurface of $\mathbf{R}^{5}_{2}$, indeed,
$V=\{X\in\mathbf{R}^{5}_{2}\mid<X,X>=-1\}$.

In the following, we  adopt  so called sausage coordinate for the
anti-de Sitter space $V$, namely,  any
$X=(x_{1},x_{2},x_{3},x_{4},x_{5})\in V$ can be expressed by

\begin{eqnarray}
\left \{\begin{array}{lll}
x_{1}&=&\frac{2r}{1-r^{2}}\sin\theta\cos\phi\\
x_{2}&=&\frac{2r}{1-r^{2}}\sin\theta\sin\phi\\
 x_{3}&=&\frac{2r}{1-r^{2}}\cos\theta\\
 x_{4}&=&\frac{1+r^{2}}{1-r^{2}}\cos t\\
x_{5}&=&\frac{1+r^{2}}{1-r^{2}}\sin t.
\end{array}\right.\nonumber\end{eqnarray}
where angular coordinates have their usual range, while $0\leq r<1$.
In this coordinates, the Lorentz metric of $V$ is
$$ds^2 = -(\frac{1+r^2}{1-r^2})^2 dt^2 + \frac{4}{(1-r^2)^2}
(dr^2 + r^2 d\theta^2 + r^2 \sin^2\theta d\phi^2),$$ thus, $t$ can
be viewed as a time function in $V$. For any slice in $V$, we can
define its height function by restriction $t$ on it. Let $M$ be a
slice of $V$, then its height function $u$ can be regarded as a
function on $\mathbf{H}^3$, which is still denoted by $u$. In the
sequel, we always assume $u$ is at least continuous at the
infinity boundary of $\mathbf {H}^3$, thus, we may define
$$w(\theta, \phi)=\lim_{r\rightarrow1}u(r, \theta, \phi),$$
hence, $w$ is a function on $\mathbf{S}^2$. Indeed, we have

\begin{theo}
Suppose $M$ is a maximal slice in $V$ , then $M$ is totally geodesic
if and only if there are constants $w_0$, $A$, $B$, $C$ on
$\mathbf{S}^2$ with $A^2 + B^2 +C^2<1$ such that

\begin{equation}
f(\theta,
\phi)=A\sin\theta\cos\phi+B\sin\theta\sin\phi+C\cos\theta,
\end{equation}
here $f=\cos(w+w_0)$.
\end{theo}

\begin{rema}
We would like to point out that $p=(\sin\theta\cos\phi,
\sin\theta\sin\phi,\\ \cos\theta)$ can be regarded as  a point on
the standard $\mathbf{S}^2\subset\mathbf{R}^3$, and each coordinate
component is a first eigenfunction of Laplacian operator on
$\mathbf{S}^2$.
\end{rema}

\begin{proof}

Suppose $M$ is a totally geodesic spacelike silce in $V$. Hence $M$
 can also  be viewed as a spacelike submanifold in $\mathbf{R}^{5}_{2}$. We
take  orthogonal frame $\{e_{1},e_{2},e_{3},e_{4},e_{5}\}$ for
$\mathbf{R}^{5}_{2}$ such that $e_{1},e_{2},e_{3},e_{4}$ and
$e_{1},e_{2},e_{3}$ are tangent vectors of  $V$ and $M$
respectively. Denote the position vector of $V$ by $X$. We may
assume  that $X=e_{5}$. Note that $M$ is totally geodesic in $V$ and
$V$ is totally umbilical in $\mathbf{R}^{5}_{2}$, we get
\begin{eqnarray}
<\widetilde{\nabla}_{e_{i}}e_{4},e_{j}>=<\widetilde{\nabla}_{e_{i}}e_{4},e_{5}>=0,
\nonumber\end{eqnarray} for $i,j=1,2,3$, here, $\widetilde{\nabla}$
is the connection in $\mathbf{R}^{5}_{2}$. Thus,  we conclude that
$e_{4}\mid_{M}=a$, where
$a=(a_{1},a_{2},a_{3},a_{4},a_{5})\in\mathbf{R}^{5}_{2}$ is a
constant vector with $<a, a>=-1$. Furthermore, one have
\begin{eqnarray}
<X\mid_{M},a>=<X\mid_{M},e_{4}>=0, \end{eqnarray} i.e., $M$ is the
intersection of $V$ and a hyperplane $\Pi_a:=\{x\in
\mathbf{R}^{5}_{2}| <x, a>=0\}$.

By (13), we obtain
\begin{eqnarray}
0=<X\mid_{M},a>=a_{1}\frac{2r}{1-r^{2}}\sin\theta\cos\phi+a_{2}\frac{2r}{1-r^{2}}\sin\theta\sin\phi
\nonumber\\+a_{3}\frac{2r}{1-r^{2}}\cos\theta-a_{4}\frac{1+r^{2}}{1-r^{2}}\cos
t-a_{5}\frac{1+r^{2}}{1-r^{2}}\sin t.\nonumber\end{eqnarray} Let
$r\rightarrow1$, we have
$$\cos(t+w_{0})=A\sin\theta\cos\phi+B\sin\theta\sin\phi+C\cos\theta,$$
or equivalently,
$$f(\theta, \phi)=A\sin\theta\cos\phi+B\sin\theta\sin\phi+C\cos\theta,$$
here,
$A=\frac{a_{1}}{\sqrt{a^{2}_{4}+a^{2}_{5}}},B=\frac{a_{2}}{\sqrt{a^{2}_{4}+a^{2}_{5}}},
C=\frac{a_{3}}{\sqrt{a^{2}_{4}+a^{2}_{5}}}$ and $\cos
w_{0}=\frac{a_{4}}{\sqrt{a^{2}_{4}+a^{2}_{5}}}$.

Conversely, if $M$ is a maximal slice, and its boundary data
satisfies (12), then we choose two constants $a_4$, $a_5$ with
$a^2_4+a^2_5
>1$ and
$$\cos w_0 = \frac{a_4}{\sqrt{a^2_4 + a^2_5}},$$
$$-\sin w_0 = \frac{a_5}{\sqrt{a^2_4 + a^2_5}}.$$
Let
\begin{equation}
\begin{split}
&a_1= \sqrt{a^2_4 + a^2_5}A\\
&a_2= \sqrt{a^2_4 + a^2_5}B\\
&a_2= \sqrt{a^2_4 + a^2_5}C.
\end{split}\nonumber
\end{equation}
Set $a=(a_1, a_2, a_3, a_4, a_5)\in\mathbf{R}^{5}_{2}$ and
$\it{C}=\{(1, \sin\theta \cos \phi, \sin\theta \sin\phi, \cos
\theta,\\ w(\theta, \phi))|0<\theta<\pi, 0\leq \phi \leq 2\pi\}$,
by a direct computation, we see that $\it{C}\subset \Pi_a$, hence
it is the boundary of $\Pi_a \cap V$ which is totally geodesic
slice of $V$, in particular, it is maximal slice, therefore, its
height function satisfies equation (1), by maximality of $M$, we
see that the height function of $M$ also satisfies the same
equation, and they are equal at the infinity boundary of
$\mathbf{H}^3$, thus, by maximal principle, we see they are equal
on $\mathbf{H}^3$ which implies $M$ is totally geodesic. Thus, we
finish to prove the theorem.
\end{proof}

As a corollary, we have
\begin{coro} Let $M$ be totally geodesic slice in $V$, then there is
a constant $w_0$ on $\mathbf{S^{2}}$ with
$$f^{2}+|\nabla^{\mathbf{S^{2}}}
f|^{2}=C,$$ where $f=\cos(w+w_0)$, $\nabla^{\mathbf{S^{2}}}$ is
the connection and  $C$ is a constant on $\mathbf{S}^2$.
\end{coro}

Combine with this fact and Theorem 3.1, we see that the Bernstein
Theorem in $V$ fails.

\bibliographystyle{amsplain}

\end{document}